\input amstex 
\documentstyle{amsppt}
\input bull-ppt
\keyedby{bull203/car}
\define\ca{\operatorname{Cap}}

\topmatter
\cvol{26}
\cvolyear{1992}
\cmonth{April}
\cyear{1992}
\cvolno{2}
\cpgs{245-252}
\ratitle
\title a general correspondence between\\ Dirichlet forms 
and right processes\endtitle
\author Sergio Albeverio and Zhi-Ming Ma\endauthor
\shorttitle{Dirichlet forms and right processes}
\address CERFIM, Cas. Post. 741, CH 6900 Locarno,
Switzerland\endaddress
\curraddr Faculty of Mathematics and SFB 237, Ruhr 
University, D 4630 Bochum,
Federal Republic of Germany          \endcurraddr
\address BiBos, University of Bielefeld, D 4800 Bielefeld, 
Federal Republic of
Germany\endaddress
\curraddr Institute of Applied Mathematics, Academia 
Sinica, Beijing 100080,
China\endcurraddr
\date May 23, 1990 and, in revised form, March 19, 
1991\enddate
\subjclassrev{Primary 31C25; Secondary 60J45, 60J25, 
60J40, 60J35}
\endtopmatter

\document
\heading 1. Introduction\endheading 

The theory of Dirichlet forms as originated by 
Beurling-Deny and developed
particularly by Fukushima and Silverstein, see e.g.\ 
\cite{Fu3, Si}, is a
natural functional analytic extension of classical (and 
axiomatic) potential
theory. Although some parts of it have abstract measure 
theoretic versions, see
e.g.\ \cite{BoH} and \cite{ABrR}, the basic  general 
construction of a Hunt
process properly associated with the form, obtained by 
Fukushima \cite{Fu2} and
Silverstein \cite{Si} (see also \cite{Fu3}), requires the 
form to be defined on
a locally compact separable space with a Radon measure $m$ 
and the form to be 
regular (in the sense of the continuous functions of 
compact support being
dense in the domain of the form, both in the supremum norm 
and in the natural
norm given by the form and the  $L^2(m)$-space). This 
setting excludes infinite
dimensional situations.

In this letter we announce that there exists an extension of
Fukushima-Silverstein's construction of  the associated 
process to the case
where the space is only supposed to be metrizable and the 
form is not required
to be regular. We shall only summarize  here results and 
techniques, for
details we refer to \cite{AM1, AM2}. Before we start 
describing our results let
us mention that some work on associating strong Markov 
processes to nonregular
Dirichlet forms had been done before, by finding a 
suitable  representation of
the given nonregular form as a regular Dirichlet form on a 
suitable
compactification of the original space. In an abstract 
general setting this was
done by Fukushima in \cite{Fu1}. The case of local 
Dirichlet forms in  infinite
dimensional  spaces, leading to associated diffusion 
processes, was studied
originally by Albeverio and H\o egh-Krohn in a rigged 
Hilbert space setting
\cite{AH1--AH3}, under a quasi-invariance and smoothness 
assumption on
$m$. This work was extended by Kusuoka  \cite{Ku} who 
worked in a Banach space
setting. Albeverio and R\"ockner \cite{AR\"O1--4} found a 
natural setting in a
Souslin space, dropping the quasi-invariance assumption. 
They also derived the
stochastic equation satisfied by the process (for 
quasi--every initial
condition) (the  ultimate result in this direction is 
contained in \cite{AR\"O5},
where also a compactification procedure by Schmuland 
\cite{Sch} and a tightness
result of Lyons-R\"ockner \cite{LR} are used). For further 
results in infinite
dimensional local Dirichlet forms see also \cite{BoH, So1, 
So2, FaR, Tak, Fu4}.
The converse program of starting from a ``good'' Markov 
process and associating
to it a, nonnecessarily regular, Dirichlet form has been 
pursued by Dynkin
\cite{D1, D2}, Fitzsimmons \cite{Fi1}, \cite{Fi2}, 
Fitzsimmons and Getoor
\cite{FG}, Fukushima \cite{Fu5}, and Bouleau-Hirsch 
\cite{BoH}. For further work
on  general Dirichlet forms see also Dellacherie-Meyer 
\cite{DM},
Kunita-Watanabe \cite{KuW}, and Knight \cite{Kn}.

Our approach differs from all the above treatments, in 
that we construct {\it
directly\/} a {\it strong Markov process\/} starting from 
a given Dirichlet
form without the assumption of regularity. In fact we 
manage to extend the
construction used in \cite{Fu3, Chapter 6}, to the 
nonregular case. By so doing
we obtain necessary and sufficient conditions for the  
existence of a certain
right process, an $m$-perfect process as explained below, 
properly associated
with the given, regular or nonregular, Dirichlet form. Our 
construction relies
on a technique we have developed in \cite{AM1} to 
associate a quasi-continuous
kernel, in a sense explained below, to a given semigroup 
(this is related to
previous work by Getoor \cite{Ge1} and Dellacherie-Meyer 
\cite{DM}). Our method of
construction of the process is related to some work by 
Kaneko \cite{Ka} who
constructed Hunt processes by using kernels which are 
quasi-continuous with
respect to a $C_{r,p}$-capacity. We also  mention that our 
work provides an
extension of a result by Y. Lejan, who obtained in 
\cite{Le1, Le2} a
characterization of semigroups associated with Hunt 
processes: in fact Lejan's
work \cite{Le1, Le2} provided an essential idea for our 
proof of the necessity
of the condition (ii) in the main theorem below.

\heading 2. Main results\endheading

We shall now present briefly our main results. Let $X$ be 
a metrizable
topological space with the $\sigma$-algebra $\scr X$ of 
Borel subsets. A
cemetery point $\Delta\notin \scr X$ is adjoined to $\scr 
X$ as  an isolated
point of $X_\Delta \equiv X\cup\{\Delta\}$. Let $(X_t 
)=(\Omega,\scr M,\scr M_t,X_t,
\Theta_t,P_x)$ be a strong Markov process with state space 
$(X_\Delta,\scr X_
\Delta)$ and life time $\zeta\equiv\inf\{t\geq 
0|X_t=\Delta\}$, where,  as in
the usual notations of e.g.\ \cite{BG}, $(\Omega,\scr M)$ 
is a measurable
space, $(\scr M_t,t\in[0,\infty])$ is an increasing family 
of sub
$\sigma$-algebras of $\scr M$, $\Theta_t$ being the shift 
and $P_x$ being the
``start measure'' (i.e.\ the measure for the paths 
conditioned to start at $x)$
on $(\Omega,\scr M)$, for each $x\in X_\Delta$. We denote 
by $(P_t)$ the
transition function of $(X_t)$ and by $(R_\alpha)$, the 
resolvent of $(X_t)$,
i.e.\ 
$$P_t f(x)=E_x[f(X_t)]$$
and
$$R_\alpha F(x)=E_x\left[\int^\infty_0e^{-\alpha 
t}f(X_t)\,dt\right],$$
where $x\in X$, $E_x$ being the expectation with respect 
to $P_x$, for all
functions $f$ for which the right-hand sides make sense.

We call $X_t$ a {\it perfect process\/} if it satisfies 
the following
properties:
\roster
\item"(i)" \<$X_t$ has the normal property: 
$P_x(X_0=x)=1$, $\forall x\in X_\Delta$;
\item"(ii)" \<$X_t$ is right 
continuous: $t\to X_t(\omega)$ is a right continuous 
function
from $[0,\infty)$ to $X_\Delta$, $P_x$ a.s., for all $x\in 
X_\Delta$;
\item"(iii)" \<$X_t$ has left limits up to 
$\zeta\:\lim_{s\uparrow t}X_s(\omega)\coloneq
X_{t-}(\omega)$ exists in $X$, for all $t\in 
(0,\zeta(\omega))$, $P_x$-a.s., $
\forall x\in X$;
\item"(iv)" \<$X_t$ has a regular resolvent in the sense 
that $R_1f(X_{t-})
I_{\{t<\zeta\}}$ is $P_x$-indistin\-guishable from 
$(R_1f(X_t))_-I_{\{t<\zeta\}}$
for all $x\in X$ and for all $f$ in the space $b\scr X$ of 
all
bounded $\scr X$-measurable functions. We have set
$(R_1f(X_t))_-I_{\{t<\zeta\}}\coloneq\lim_{s\uparrow
t}R_1f(X_s)I_{\{t<\zeta\}}$ (where we always make the 
convention that
$Z_{0-}=Z_0$ for any process $Z_t$, $t\geq 0)$.\endroster

\rem{Remarks} (i) A strong Markov process satisfying only 
(i), (ii) will be
called a {\it right process\/} with Borel transition 
semigroup. If $X$ is a
Radon space this definition coincides with the one in 
\cite{Sh, Definition
(8.1)} and \cite{Ge2, (9.7)}.

(ii) A special standard process in the sense of \cite{Ge2, 
(9.10)} and, in
particular, a Hunt process (cf.\ e.g.\ \cite{BG}), always 
satisfy (iii), (iv),
hence they are particular perfect processes (see \cite{Ge2 
(7.2), (7.3)},
\cite{Sh (4.7.6) and (9.7.10)} for the proof).

Thus we have the following inclusions, cf.\ \cite{Ge2, p.\ 
55}:
$$\text{(Feller)}\subset\text{(Hunt)}\subset\text{(special 
standard)}\subset
\text{(perfect)}\subset\text{(right)}\.$$
\endrem

In what follows we shall assume that $m$ is a 
$\sigma$-finite Borel measure on
$X$. A process $(X_t)$ is called {\it $m$-tight\/} if 
there exists an
increasing sequence of compact sets $\{K_n\}$, $n\in \bold 
N$ of $X$ such that
$$P_x\{\lim_n\sigma_{X-K_n}\geq \zeta\}=1,\qquad 
m\text{-a.e.\ }x\in X,$$
where for any subset $A$ of $X_\Delta$, $\sigma_A\coloneq 
\inf\{t>0\:X_t
\in A\}$ is the hitting time of $A$.

A process $(X_t)$ is called an {\it $m$-perfect process\/} 
if it is a perfect
process and is $m$-tight. It follows from an idea of T. J. 
Lyons and M.
R\"ockner \cite{LR} that any strong Markov process $(X_t)$ 
in a metrizable Lusin space is
$m$-tight if it satisfies (ii), (iii), see \cite{AMR}. 
Thus in a metrizable Lusin
the concepts of perfect process and $m$-perfect process 
coincide. (In the special case
of $(X_t)$ being a standard process on a  locally compact 
metrizable space the
above conclusion that $(X_t)$ is $m$-tight can also be 
derived from \cite{BG,
(9.3)}.)

Let us now give the correlates of above definitions for 
Dirichlet forms. Let 
$(\scr E ,\scr F )$ be a Dirichlet form on $L^2(X,m)$, 
i.e. $\scr E $ is a
positive, symmetric, closed bilinear form on $L^2(X,m)$ 
such that unit
contractions operate on $\scr E $ (i.e.\ $\scr E 
(f^{\#},f^{\#})\leq \scr E
(f,f)$, $f^{\#}=(f\vee 0)\wedge 1$, all $f$ in the 
definition domain $\scr F $
of $\scr E )$. We set as usual $\scr E _1(f,g)\equiv \scr 
E (f,g)+(f,g)$, $
\forall f$, $g\in \scr F $, where $(f,g)$ is the 
$L^2(X,m)$-scalar product of
$f$ and $g$.

In the sequel we always regard $\scr F $ as a Hilbert 
space equipped with the
inner product $\scr E _1$. For any closed subset $F\subset 
X$, we set
$$\scr F _F\equiv \{f\in F\:f=0\ m\text{-a.e.\ on }X-F\},$$
$\scr F _F$ is then a closed subset of $\scr F $.

The following definitions are extensions of corresponding 
definitions in
\cite{Fu3}. An increasing sequence of closed sets 
$\{F_k\}_{k\geq 1}$ is called
an {\it $\scr E $-nest\/} if $\bigcup_k\scr F _{F_k}$ is 
$\scr E _1$-dense in
$\scr F $. A subset $B\subset X$ is said to be {\it $\scr 
E$-polar\/} if there
exists an $\scr E$-nest $\{F_k\}$ such
that
$$B\subset \bigcap_k(X-F_k)\.$$
A function $f$ on $X$ is said to be {\it $\scr E 
$-quasi-continuous\/} if there
exists an $\scr E $-nest $\{F_k\}$ such that $f_{|F_k}$, 
the restriction of $f$
to $F_k$, is continuous on $F_k$ for each $k$.

It is not difficult to show that every $\scr E $-polar set 
is $m$-negligible.
We denote by $T_t$ resp.\ $G_\alpha$ the semigroup resp.\ 
resolvent on
$
L^2(X,m)$ associated with $(\scr E,\scr F)$; i.e., if $A$ 
is the  generator of
$T_t$,
$$\left(\sqrt{-A} f, \sqrt{-A} g\right)=\scr E (f,g)\quad 
\forall f,\ g\in \scr
F =D\left(\sqrt{-A}\right)\.$$
(cf.\ \cite{Fu3}). We set
$$\scr H=\{h\:h=G_1f\text{ with }f\in L^2(X;m),0<f\leq 
1m\text{-a.e.}\}$$
We remark that $\scr H$ is nonempty because we assumed $m$ 
to be
$\sigma$-finite.  For $h\in \scr H$ we define the {\it 
$h$-weighted capacity\/}
$\ca_h$ as follows
$$
\ca_h(G)\coloneq\inf\{\scr E _1(f,f)\:f\in \scr F ,f\geq 
hm\text{-a.e.\ on
}G\},$$
for any open subset $G$ of $X$, and
$$\ca_h(B)\coloneq \inf\{\ca_h(G)\:G\supset B,G\text{ 
open}\}$$
for any arbitrary subset $B\subset X$. It is possible to 
show, see \cite{AM2},
that $\ca_h$ is a Choquet capacity enjoying the important 
property of countable
subadditivity. The relation between  this notion of 
capacity and the notion of
$\scr E $-nest is expressed by the following proposition:
\proclaim{Proposition} An increasing sequence of closed 
subsets $\{F_k\}$ of
$X$  is an $\scr E $-nest if and only if for some $h\in 
\scr H$ \RM(hence for
all $h\in \scr H)$ one has
$$\ca_h(X\bsl F_k)\downarrow0\quad \text{ as }k\to\infty\.$$
\endproclaim
\noindent For the proof we refer to Proposition 2.5 of 
\cite{AM2}.

Denote by Cap the usual 1-capacity given by $\scr E $, 
cf.\ \cite{Fu3}. One has
obviously  $\ca_h(B)\leq \ca(B)$ for every $B\subset X$. 
Hence the following
corollary holds:
\proclaim{Corollary} Every set $B\subset X$ with $\ca$ 
$B=0$ is an $\scr E
$-polar set. Every nest $\{F_k\}$ resp.\ every 
quasi-continuous function in the
sense of\ \ \cite{Fu3} is an $\scr E $-nest resp.\ an 
$\scr E $-quasi-continuous
function in  our sense \RM(we remark that \cite{Fu3} the 
space $X$
is supposed to be locally compact separable and $m$ to be 
supported by $X$ and
Radon\RM).\endproclaim

Now let $(X_t)$ be a Markov process with transition 
function $P_t$. We say that
$(X_t)$ is {\it associated with\/} $\scr E $ if 
$P_tf=T_tf$ $m$-a.e.\ for all
$f\in L^2(X,m)$, $t>0$, and it is properly {\it associated 
with\/} $\scr E$ if $P_tf$ is an $\scr E$-quasi-continuous
version of $T_tf$ for all $f\in L^2(X,m)$, $t>0$. The main
result we obtain is the following
\proclaim{Theorem} Let $(\scr E,\scr F)$ be a Dirichlet 
form on $L^2(X;m)$.
Then the following family of conditions \RM{(i)--(iii)} is 
necessary and
sufficient for the existence of an $m$-perfect process 
$(X_t)$
associated with $\scr E $\RM:
\roster
\item"(i)" there exists an $\scr E $-nest $\{X_k\}$ 
consisting of compact subsets
of $X$\RM;
\item"{(ii)}" there exists an $\scr E _1$-dense subset 
$\scr F _0$
of $\scr F $ consisting of $\scr E $-quasi-\linebreak 
continuous functions\RM;
\item"{(iii)}" there exists a countable subset $B_0$ of 
$\scr F _0$ and an 
$\scr E$-polar subset $\scr N$ such that
$$\sigma\{u\:u\in B_0\}\supset \scr 
B(X)\cap(X-N)\.$$\endroster

Moreover, if an $m$-perfect process $(X_t)$ is associated 
with $\scr E $, then
it is always properly associated with $\scr E $.\endproclaim
\rem{Remarks} (i) If $\scr E $ is a regular  Dirichlet 
form in the sense of
\cite{Fu3} then all conditions are satisfied. But the 
regularity assumption on
the Dirichlet form $\scr E $ usually assumed in the 
literature, cf.\ \cite{Fu3,
Si}, is not necessary for the existence of  an $m$-perfect 
process
(it is even not necessary for the existence of a diffusion 
process). Assume in
fact each single point set of $X$ is a set of zero 
capacity (e.g.\
 $X= \bold R^d$, $d\geq 2$, $\scr E $ the classical 
Dirichlet form associated
with the Laplacian on $\bold R^d)$. Let $\mu$ be a smooth 
measure (in the sense 
of \cite{Fu3}), which is nowhere Radon (i.e.\ 
$\mu(G)=\infty$ for all nonempty
open subsets $G\subset X)$: the existence of such nowhere 
Radon smooth measures
has been proven in \cite{AM4} (see also \cite{AM3, AM5}).

We consider the perturbed form $(\scr E ^\mu,\scr F ^\mu)$ 
defined as follows:
$$\scr F ^\mu\coloneq \scr F \cap L^2(X;m);\qquad \scr E 
^\mu(f,g)
\coloneq \scr E(f,g)+\int_Xfg\mu(dx)\quad 
\forall f,\ g\in \scr F ^\mu\.$$
It has been proven in \cite{AM5, Proposition 3.1} that 
$(\scr E^\mu,\scr F ^\mu)$ is
again a Dirichlet form. One can check that $(\scr E ^\mu 
,\scr F ^\mu)$
satisfies all conditions in the theorem, see \cite{AM7}, 
hence the theorem is
applicable and there exists an $m$-perfect process 
associated with 
$(\scr E ^\mu ,\scr F ^\mu)$. Moreover, if $(\scr E,\scr 
F)$ satisfies the
local property, then so does $(\scr E ^\mu ,\scr F ^\mu)$ 
and hence there
exists a diffusion process associated with $(\scr E ^\mu 
,\scr F ^\mu)$ (see
(ii) below). But clearly $(\scr E ^\mu ,\scr F ^\mu)$ is 
not regular, in fact
there is even no nontrivial continuous function belonging 
to $\scr F ^\mu$,
since $\mu$ is nowhere Radon. See \cite{AM7} for details.

(ii) An $m$-perfect process is a diffusion (i.e.\ 
$P_x\{X_t$ is continuous in
$t\in[0,\zeta)\}=1$, for q.e.\ $x\in X\})$ if and only if 
the associated
Dirichlet form $(\scr E,\scr F)$ satisfies the local 
property (in the sense
that $\scr E (f,g)=0$ if $\supp f\cap \supp g=\emptyset)$, 
see \cite{AM6}.

(iii) By requiring $\scr F _0$ in (ii) to consist of 
strictly $\scr E
$-quasi-contin\-uous functions we obtain a  necessary and 
sufficient condition
for the existence of a Hunt process associated with $(\scr 
E,\scr F)$, see
\cite{AM8}.

(iv) By introducing a dual $h$-weighted capacity and 
employing a Ray-Knight
compactification method it is possible to extend the above 
theorem to
nonsymmetric Dirichlet forms satisfying the sector 
condition.

(v) Applications of the above theorem to infinite 
dimensional spaces $X$ are in
preparation. They allow in particular to construct 
infinite dimensional
processes with discontinuous sample paths, with 
applications to certain systems
with infinitely many degrees of freedom.
\endrem
\heading Acknowledgments\endheading

We are very indebted to Masatoshi Fukushima, who greatly 
encouraged our work,
and carefully read the first version of the manuscripts of 
\cite{AM1--2}
suggesting many improvements. We are also very grateful to 
H. Airault,
J. Brasche, P. J. Fitzsimmons, R. K. Getoor, W. Hansen, M. 
L.
Silverstein, S. Watanabe, R. Williams, J. A. Yan, Zhang, and
especially Michael R\"ockner for very interesting and 
stimulating
discussions. We also profited from meetings in Braga and 
Oberwolfach and are
grateful to M. De Faria, L. Streit resp., H. Bauer, and  
M. Fukushima for
kind  invitations. Hospitality and/or financial support by 
BiBoS, A. von
Humboldt Stiftung, DFG, and Chinese National Foundation is 
also gratefully
acknowledged.

\subheading\nofrills{{\it Note added in Proof}}. After we 
finished this work we learned
from P. J. Fitzsimmons that every (nearly) $m$-symmetric 
right process is an
$m$-special standard process (see \cite{Fi1}), and our 
notion of
perfect process is equivalent to a special standard 
process. Thus
our theorem gives in fact a general correspondence between
Dirichlet forms and right processes. We are most grateful 
to P. J. Fitzsimmons
for his remark.

\enddocument